\documentclass[12pt]{amsart}
\usepackage{amsmath, amssymb, amsthm}

\newcommand{\zn}{\mathbb{Z}_n}

\DeclareMathOperator{\soc}{Soc}

\newtheorem{theorem}{Theorem}[section]

\newtheorem{lemma}[theorem]{Lemma}

\theoremstyle{definition}
\newtheorem{definition}[theorem]{Definition}
\newtheorem{example}[theorem]{Example}

\theoremstyle{remark}
\newtheorem{remark}[theorem]{Remark}

\numberwithin{equation}{section}

\begin{document}

\title{Test ideals in diagonal hypersurface rings II}
\author{Moira A. McDermott}
\address{Department of Mathematics and Computer Science,
Gustavus Adolphus College,
Saint Peter, MN 56082}
\email{mmcdermo@gac.edu}
\date{\today}

\begin{abstract}
Let $R=k[x_1, \ldots, x_n]/(x_1^d + \cdots + x_n^d)$, where $k$ is a field
of characteristic $p$, $p$ does not divide $d$ and $n \geq 3$.
We describe a method for computing the test ideal for these diagonal
hypersurface rings.
This method involves using a characterization of test ideals in
Gorenstein rings as well as developing a way to compute tight closures
of certain ideals despite the lack of a general algorithm.
In addition, we compute examples of test ideals in diagonal
hypersurface rings of small characteristic (relative to $d$) including
several that are not integrally closed.  These examples provide a
negative answer to Smith's question of whether the test ideal in
general is always integrally closed \cite{multideal}.
\end{abstract}

\maketitle

\section*{Introduction}
Test ideals play an important role in the theory of tight closure developed
by Melvin Hochster and Craig Huneke.  
Unfortunately, both tight closure and test ideals are difficult to
compute in general.
In this paper we describe a method for computing the test ideal for
diagonal hypersurfaces
$k[x_{1}, \ldots , x_{n}]/(x_{1}^{d} + \cdots + x_{n}^{d})$, where $k$ 
is a field of characteristic $p$, $p$ does not divide $d$ and $n \geq 3$.
This method involves using a characterization of test
ideals in Gorenstein rings as well as developing a way to compute
tight closures of certain ideals despite the lack of a general algorithm.

It is worth noting that the test ideals for these diagonal hypersurface rings
are very different depending on the magnitude of $p$ (usually relative
to $d$).
If $d\geq n$ and the dimension of the ring is two, then results of
Huneke and Smith \cite{HS} show that the test ideal is 
$(x_1, x_2, x_3)^{d-3+1}$ for $p>d$.
In this case the test ideal is essentially all elements of degree
greater than the a-invariant.  This follows from their tight closure
interpretation of the Strong Kodaira Vanishing Theorem.  Huneke and
Smith also point out that the Vanishing Theorem is true for hypersurface rings 
with a similar restriction on the characteristic.
Huneke gives a direct proof of the Strong Vanishing Theorem for hypersurfaces
\cite[6.4]{Hu2} using ideas found in earlier work of Fedder \cite{Fe87}.  
There are also similar results in \cite[Cor 3]{Ha3} using the idea of
F-injectivity in negative degree.
For the diagonal hypersurfaces, the Strong Vanishing Theorem implies that
if $d\geq n$, then the test ideal is
$(x_1, \ldots, x_n)^{d-n+1}$ for sufficiently large $p$, and if $d<n$,
then the test ideal is the unit ideal \cite[6.3]{Hu2}.
Fedder and Watanabe also have results that show that if $d<n$, then
the ring is F-regular \cite[(2.11)]{FW} and hence the test ideal is  the
unit ideal, again for sufficiently large $p$.
On the other hand, we show in \cite{testideal1} that if $p<d$, then the test
ideal is contained in $(x_1, \ldots, x_n)^{p-1}$, which is much
smaller than would be expected in many cases.  This result does not
depend on $n$.  We also show in \cite{testideal1} that if $p=d-1$,
then the test ideal is in fact equal to $(x_1, \ldots, x_n)^{p-1}$.

We are interested in computing test ideals in these rings
when $p$ is less than the bounds required in \cite{Fe87}, \cite{Hu2},
\cite{Ha3} and \cite{FW}.
For $k[x_{1}, \ldots , x_{n}]/(x_{1}^{d} + \cdots + x_{n}^{d})$, where $k$ 
is a field of characteristic $p$, the bound in \cite{Fe87} is 
$p>n(d-1)-d$.  The bound in \cite{Hu2} is $p>n(d-1)-2d+1$ and the bound
in \cite{Ha3} is $p>n(d-1)-2d$.  The bound in \cite{FW} comes from the
bound in \cite{Fe87} and so is also $p>n(d-1)-d$.  
It is quite likely that these are not the best possible bounds.
The bound in the two dimensional case, $p>d$ \cite{HS}, is quite a bit
better than these.  In many examples when the dimension
is greater than two, the bound of $p>d$ is sufficient.
For this reason, and the fact that the case when $p=d-1$ is known
\cite{testideal1}, we are particularly interested in computing test
ideals when $p<d-1$.
We do have one example (\ref{n=5,d=4,p=7}) where $p$ is greater than $d$
but less than the bound in \cite{FW} ($p>n(d-1)-d$), and the ring is not
$F$-regular as predicted.

In this paper we describe a method for computing test ideals in
diagonal hypersurface rings.  We then use this method to compute many
examples of test ideals when $p<d-1$ and when $p$ less than the
previously mentioned bounds.

Recently Karen Smith has shown that
test ideals are closely related to certain multiplier ideals that arise in
vanishing theorems in algebraic geometry.  
In \cite{multideal} Smith established that for 
normal local Cohen Macaulay $\mathbb Q$-Gorenstein rings essentially
of finite type over a field of characteristic zero,
the multiplier ideal 
is a {\it universal test ideal} 
(see \cite {Na} and \cite{Ein} for more information about the
multiplier idea).  
Similar results were also obtained independently by Hara in \cite{Ha2}.
The multiplier ideal is integrally
closed.  This lead Smith to ask whether the test ideal in general is
integrally closed.  Several of the examples we are able to compute are
particularly interesting because they are not integrally closed.
These examples provide a partial negative answer to Smith's question.

\section{Notation and Definitions}

Throughout this paper $R$ is a commutative Noetherian ring of
prime characteristic $p>0$.  The letter $q$ will always stand for a
power $p^e$ of $p$, where $e \in \mathbb N$.

We review the definition of tight closure for ideals of rings of 
characteristic $p > 0$.  Tight closure is defined more generally for 
modules and also for rings containing fields of arbitrary 
characteristic.  See \cite{HH1} or \cite{Hu1} for more details.

\begin{definition}\label{deftc}
Let $R$ be a ring of characteristic $p$  and 
$I$ be an ideal in a Noetherian ring $R$ of characteristic $p > 
0$.  An element $u \in R$ is in the {\it tight closure} of $I$, 
denoted $I^{*}$,
if there exists an element 
$c \in R$, not in any minimal prime of $R$, such that for all large 
$q=p^{e}$, $cx^q \in I^{[q]}$ 
where $I^{[q]}$ is the ideal generated by the $q$th powers of all 
elements of $I$.
\end{definition}

In many applications one would like to be able to choose the element
$c$ in the definition of tight closure independent of $x$ or $I$.  It
is very useful when a single choice of $c$, a test element, can be
used for all tight closure tests in a given ring.

\begin{definition}\label{deftestideal}
The ideal of all $c \in R$ such that for any ideal $I \subseteq R$, 
we have $cu^q \in I^{[q]}$ for all $q$ whenever
$u \in I^*$ is called the {\it test ideal} for $R$ and is denoted by $\tau$.  
An element of the test ideal that is not in any minimal prime is 
called a {\it test element}.
\end{definition}

\section{Determining the Test Ideal}

We will make use of the following grading in our calculation of the 
test ideal.
We denote by $\mathbb{Z}_n$ the ring $\mathbb{Z} / n\mathbb{Z}$.
Next we describe a $\mathbb{Z}_n$-grading of rings of the form 
$R = A[z]/(z^n - a)$ where $a \in A$.  The ring $R$ has the following
decomposition as an $A$-module:
$$
R = A \oplus Az \oplus \cdots \oplus Az^{n-1}.
$$
This is true because every element of $R$ can be uniquely expressed as
an element of $A \oplus Az \oplus \cdots \oplus Az^{n-1}$ by replacing every
occurrence of $z^n$ by $a$.  $R$ is $\mathbb{Z}_n$-graded, where the $i$th
piece of $R$, denoted by $R_{i}$, is $Az^i$, $0 \leq i < n$, since 
$Az^iAz^j \subseteq Az^{i+j}$ if $i+j <n$ 
and 
$Az^iAz^j \subseteq Az^{i+j-n}$ if $i+j \geq n$.

We use this idea to obtain multiple $\mathbb{Z}_n$-gradings of 
$R = k[x_{1}, \ldots , x_{n}]/
$\linebreak$
(x_{1}^d + \cdots + x_{n}^d)$, 
where $k$ is a field of characteristic $p$.
Let $z=x_{i}$ and 
$A=k[x_{1}, \ldots , \hat{x_{i}}, \ldots , x_{n}]$.

It is not difficult to show that if $I$ is a graded ideal, then so is $I^*$.

\begin{lemma}\label{gradinglemma}\cite{testideal1}
Let $R$ be a finitely generated $k$-algebra that is 
$\mathbb{Z}_{n}$-graded and of characteristic $p$, where $p$ is not a 
prime factor of $n$ ($p = 0$ is allowed).  
Then the tight closure of a homogeneous ideal of $R$ is homogeneous.
\end{lemma}

We will use the following result about the test ideal in a Gorenstein
ring to compute the test ideal for the diagonal hypersurfaces.

\begin{lemma}\label{colon}\cite[2.14]{Hu1},\cite[4.3]{Hu2}
Let $(R,m)$ be a Gorenstein local ring with an $m$-primary test
ideal.  Let $I$ be generated by a system of parameters consisting of
test elements.  Then the test ideal is $I \colon I^*$.
\end{lemma}

We also know that the elements of the Jacobian ideal are always test 
elements  by Hochster and Huneke's ``test elements
via Lipman Sathaye'' result \cite[(1.5.5)]{HH4}, \cite[3.12]{Hu2}.
In practice, it is preferable to use
$x_{1}^{d}, \ldots , x_{n-1}^{d}$
as the sequence of test elements.
This allows us to capitalize on certain symmetries that arise in the
diagonal hypersurface rings.  Because we have the defining relation
$x_1^d + \cdots + x_n^d =0$, the ideals generated by any $n-1$ of the
$d$th powers of the variables are the same.

Let $R=k[x_{1} \ldots x_{n}]/(x_{1}^d + \cdots + 
x_{n}^d)$ where 
$k$ is a field of characteristic $p$, $p$ does not divide $d$, and 
$n \geq 3$.
Then 
$\tau = (x_1^d, \ldots, x_{n-1}^d):(x_1^d, \ldots, x_{n-1}^d)^*$ and the
problem of computing the test ideal reduces to determining
$(x_1^d, \ldots, x_{n-1}^d)^*$.

\begin{remark}
It is also worth noting that
the test ideal for $R$ can be generated by 
monomials.  Recall that there is a $\zn$-grading of $R$ associated with 
each $x_{i}$, $1 \leq i \leq n$.  We also know that if $I$ is a 
homogeneous ideal, then so are $I^{*}$ (Lemma \ref{gradinglemma}) and 
$I \colon I^{*}$.  
Since
$(x_{1}^d, \ldots , x_{n-1}^d)$ 
is homogeneous with respect to each of the gradings, so is
$(x_{1}^d, \ldots , x_{n-1}^d)
{\colon}_{R}
(x_{1}^d, \ldots , x_{n-1}^d)^{*}$ 
and hence so is the test ideal.
Using the grading with respect to each $x_{i}$, the multigrading, we 
see that the test ideal can be generated by monomials.  
Essentially, this is because only monomials are homogeneous with 
respect to all $n$ gradings simultaneously.
\end{remark}

\section{Computing $(x_1^d, \ldots, x_{n-1}^d)^*$}
Let $R=k[x_{1} \ldots x_{n}]/(x_{1}^d + \cdots + 
x_{n}^d)$ where 
$k$ is a field of characteristic $p$, $p$ does not divide $d$, and 
$n \geq 3$.
First we describe a general method for computing 
$(x_1^d, \ldots, x_{n-1}^d)^*$.
Let
$m=(x_1, \ldots, x_n)$ and $I=(x_1^d, \ldots, x_{n-1}^d)$.

\begin{enumerate}
\item
Let $J$ be a candidate for $I^*$.
(Begin with $J=I$.)
\item
Compute $\soc(R/J) = J\colon_R m = (u_1, \ldots, u_m)$.
\item
Determine whether $u_i \in I^*$, $1 \leq i \leq m$.
\item
Form a new candidate for $I^*$ by adding all $u_i$'s that are in $I^*$
to $J$ and repeat.
\end{enumerate}

When no generators of $\soc(R/J)$ are in $I^*$, the process is
complete and $J=I^*$.  Since $R$ is Noetherian, the process must
eventually end.  Next we explain why it is sufficient to check
elements of $\soc(R/J)$.

\begin{lemma}\label{checksocle}
Let $(R,m,k)$ be a Noetherian local ring and $J \subseteq I$ ideals of
$R$.  If $J \subset I$, then $I$ contains an element of $\soc(R/J)$.
\end{lemma}

\begin{proof}
Let $u \in I \backslash J$.  We know that $\soc(R/J)$ contains all
simple submodules of $R/J$ and therefore meets every submodule of
$R/J$.
Consider $N = \frac{Ru}{J}$, the submodule generated by $u$.  If $N
\neq (0)$, then $N \cap \soc(R/J) \neq (0)$.  Let $\overline u$ be the
image of $u$ in $R/J$.  There exists $r \in R$ with $r\overline u \in
\soc(R/J)$.  Then $ru$ is the desired element of $I$.
\end{proof}

Suppose $J$ is a candidate for $I^*$.
By Lemma \ref{checksocle}, if $I^*$ strictly contains $J$, then $I^*$ must contain an
element of $\soc(R/J)$.
Therefore, to show $J=I^*$, it is
sufficient to show that $\sum ku_i \cap I^*=0$ where $\soc(R/J)$ is generated 
by the $u_i$.

The final step in the ``algorithm'' to be justified is step 3.  The
method described above is not a true algorithm since there is no
known algorithm for computing tight closure, except in some special cases.
{\it A priori}, one might have to test infinitely many exponents in
order to determine whether one element is in the tight closure of a
given ideal.  Another way to describe this problem is in terms of test
exponents as in \cite{HHtestexp}.
 
\begin{definition}
Let $R$ be a reduced Noetherian ring of positive prime characteristic
$p$.  Let $c$ be a fixed test element for $R$.  We shall say that
$q=p^e$ is a {\it test exponent for $c$, $I$, $R$} if whenever $cu^Q
\in I^{[Q]}$ and $Q \geq q$, then $u \in I^*$.
\end{definition}

Whenever one can compute what the test exponent is, one obtains an
effective test for tight closure.  Despite the lack of test exponents,
in practice, one can compute tight closure in diagonal hypersurface
rings in many cases.
The following two observations are helpful in determining whether a
specific element is in the tight closure of a given ideal.
Even though in general there is no bound on the power of an element
needed to test tight closure, there are two situations where one
exponent is enough. 

\begin{remark}\label{frobenius}
Note that the Frobenius closure $I^F$ of an ideal $I$ is the set of
elements $u$ such that $u^q \in I^{[q]}$ for $q \gg 0$.  It follows
from elementary properties of the Frobenius map that
if $u^{q'} \in I^{[q']}$ for one value of $q'=p^e$, then $u^{q} \in I^{[q]}$
for all higher powers $q \geq q'$.  In fact, $I^F$ is often defined to
be the set of elements $u$ such that $u^q \in I^{[q]}$ for some $q=p^e$.
This means that if
$u^{q} \in I^{[q]}$ for just {\bf one} value of $q$, then
$u \in I^F \subseteq I^*$.
\end{remark}

\begin{remark}\label{notintc}
Recall that
if $c$ is a test element, then
$cu^{q} \in I^{[q]}$ for all $u \in I^*$ and for {\bf all} $q=p^e$.
This means that if
$cu^{q} \notin I^{[q]}$ for even {\bf one} $q$, then $u \notin I^*$.
\end{remark}

In principle, this method would only work if $I^F=I^*$.  There is some
evidence that $I^F=I^*$ when $n=d=3$ and $p \equiv 2 \mod 3$
\cite{TAMS}, however we are not conjecturing that $I^F=I^*$ for all
diagonal hypersurface rings.  In practice, however, the potential gap
between $I^F$ and $I^*$ has never prevented us from computing a test
ideal.  Instead, the current limitations are memory and monomial
bounds in Macaulay~2 \cite{M2}.  In almost every example, we show that an
element is in the tight closure of an ideal by using the observation
in (\ref{frobenius}) and showing that the element is actually in the
Frobenius closure of the ideal.  In a few cases, we have used the
following result of Hara.  Smith has a similar result \cite[Lemma 3.2]{Smgraded}.

\begin{lemma}\cite[Lemma 2]{Ha3}
Let $R=\bigoplus_{n\geq0}R_n$ be a Noetherian $\mathbf N$-graded ring
defined over a perfect field $k=R_0$ of characteristic $p>0$.  Assume
that $R$ is Cohen-Macaulay.  Let $x_1, \ldots, x_d$ be a homogeneous
system of parameters of $R$, and assume that $d=\dim R \geq 1$.  If a
homogeneous element $z$ satisfies $\deg(z) \geq \sum_{i=1}^d
\deg(x_i)$, then $z \in (x_1, \ldots, x_d)^*$.
\end{lemma}

It is interesting to note that every
instance where we could not show that an element was in the Frobenius
closure of an ideal by direct computation was an instance where Hara's
lemma applied.  

\section{Examples}

\begin{example}\label{n=d=5,p=2}
Let $R$ be the localization at $(x_1, \ldots, x_5)$ of the ring
\[
\frac{k[x_1, \ldots, x_5]}{(x_1^5+x_2^5+x_3^5+x_4^5+x_5^5)}
\]
where $k$ is a field of characteristic two.  In this case the test ideal 
for $R$ is generated by the elements $x_i^2x_j$ 
for all $1 \leq i,j \leq 5$.

To verify this, we use the observation that in a Gorenstein ring with 
isolated singularity, the test ideal is $J:J^*$ where $J$ is an ideal
generated by a system of parameters that are test elements (Lemma \ref{colon}).
We also know that elements in the Jacobian ideal are test elements 
\cite[(1.5.5)]{HH4}, so in this ring $x_1^4, \ldots, x_5^4$ are test elements.
Some of our calculations will be easier if we use the fact that
$x_1^5, \ldots, x_5^5$ are also test elements.
Thus we use $J=(x_1^5, \ldots, x_4^5)$.  One can calculate directly that
\begin{multline*}
J^* = (x_1^5, x_2^5, x_3^5, x_4^5, 
x_1^3x_2^3x_3^3x_4^3x_5^3,
x_1^4x_2^4x_3^4x_4^4x_5^2,
x_1^4x_2^4x_3^4x_4^2x_5^4,\\ 
x_1^4x_2^4x_3^2x_4^4x_5^4,
x_1^4x_2^2x_3^4x_4^4x_5^4,
x_1^2x_2^4x_3^4x_4^4x_5^4).\\
\end{multline*}
Using the $\mathbb{Z}_5$ grading we can assume that $J^*$ is generated by 
elements of the form $u = x_1^{a_1}x_2^{a_2}x_3^{a_3}x_4^{a_4}x_5^{a_5}$.
Also any monomial of that form in $J^*\backslash J$ must have all
$0<a_i<5$.  Clearly, we must have all $a_i<5$ in order to have $u
\notin J$.   The fact that all $a_i>0$ follows from ``tight closure
from contractions'' \cite[1.7]{Hu1} since 
$k[x_1, \ldots, x_5]/(x_1^5+\cdots+x_5^5)$
is a module finite extension of $k[x_1, \ldots , x_4]$.
To verify that the monomials listed above are in
$J^*$ we use the observation in (\ref{frobenius}), namely
if $u^{q'} \in I^{[q']}$, then $u^{q} \in I^{[q]}$, $q \geq 
q'$ and hence $u \in I^{*}$.
One easily checks that $u^4 \in J^{[4]}$ for all monomials $u$ listed above.
It is also easy to check that no generators of the socle modulo the
candidate for $J^*$ are in $J^*$.  Since $J^*$ and $m$ are both
monomial ideals, it is routine to compute $J^*:m$, the socle modulo $J^*$.
We compute $J^*:m$ and see that the socle modulo $J^*$ has 25 generators
that are not in $J^*$.  Those generators are as follows:
\begin{eqnarray*}
 &&u_1=x_1  x_2^4x_3^4x_4^4x_5^4, u_2=x_1^4x_2x_3^4x_4^4x_5^4,\ldots,
   u_5=x_1^4x_2^4x_3^4x_4^4x_5, \\
 &&u_6=x_1^2x_2^3x_3^4x_4^4x_5^4, u_7=x_1^2x_2^4x_3^3x_4^4x_5^4, \ldots,
u_{30}=x_1^4x_2^4x_3^4x_4^3x_5^2.
\end{eqnarray*}
We use $c=x_1^4$ as a test element and see that, for example,
$cu_1^{32} \notin J^{[32]}$ and $cu_6^{16} \notin J^{[16]}$.  Similar
calculations and the observation in (\ref{notintc}) show that the
remaining monomials are not in $J^*$.  Computing $J:J^*$ gives the desired result.
\end{example}

\begin{example}\label{n=4,d=7,p=3}
Let $R$ be the localization at $(x_1, \ldots, x_4)$ of the ring
\[
\frac{k[x_1, \ldots, x_4]}{(x_1^7+x_2^7+x_3^7+x_4^7)}
\]
where $k$ is a field of characteristic three.  In this case the test ideal 
for $R$ is generated by the elements 
$x_i^2x_j^2$  for all $1 \leq i,j \leq 4$.

In this example we let $J=(x_1^7, x_2^7 , x_3^7)$.
One can calculate directly that
\[
J^* = (x_1^7, x_2^7, x_3^7,
x_1^3x_2^5x_3^5x_4^5,
x_1^5x_2^3x_3^5x_4^5,
x_1^5x_2^5x_3^3x_4^5,
x_1^5x_2^5x_3^5x_4^3
).
\]
As in the previous example, we use the observation in (\ref{frobenius})
and check that $u^3 \in J^{[3]}$ for all monomials $u$ listed above.
We compute $J^*:m$ and see that the socle modulo $J^*$ has 10 generators
that are not in $J^*$.  Those generators are as follows:
\begin{eqnarray*}
&&u_1=x_1^2x_2^6x_3^6x_4^6, u_2=x_1^6x_2^2x_3^6x_4^6, 
u_2=x_1^6x_2^6x_3^2x_4^6, u_4=x_1^6x_2^6x_3^6x_4^2, \\
&&u_5=x_1^4x_2^4x_3^6x_4^6, u_6=x_1^4x_2^6x_3^4x_4^6,\ldots,
u_{10}=x_1^6x_2^6x_3^4x_4^4.
\end{eqnarray*}
We use $c=x_1^6$ as a test element and see that, for example,
$cu_1^{9} \notin J^{[9]}$ and $cu_5^{9} \notin J^{[9]}$.  Similar
calculations and the observation in (\ref{notintc}) show that the
remaining monomials are not in $J^*$.  Computing $J:J^*$ gives the desired result.
\end{example}

\begin{example}\label{n=5,d=4,p=7}
Let $R$ be the localization at $(x_1, \ldots, x_5)$ of the ring
\[
\frac{k[x_1, \ldots, x_5]}{(x_1^4+x_2^4+x_3^4+x_4^4+x_5^4)}
\]
where $k$ is a field of characteristic seven.  
In this case the test ideal for $R$ is $(x_1, \ldots, x_5)$, the
maximal ideal.  

In this example we let $J=(x_1^4, x_2^4 , x_3^4, x_4^4)$.
One can calculate directly that
\[
J^* = (x_1^4, x_2^4, x_3^4,
x_1^3x_2^3x_3^3x_4^3x_5^3
).
\]
As in the previous example, we use the observation in (\ref{frobenius})
and check that $(x_1^3x_2^3x_3^3x_4^3x_5^3)^7 \in J^{[7]}$.
We compute $J^*:m$ and see that the socle modulo $J^*$ has 5 generators
that are not in $J^*$.  Those generators are as follows:
\begin{eqnarray*}
&&u_1=x_1^2x_2^3x_3^3x_4^3x_5^3,
u_2=x_1^3x_2^2x_3^3x_4^3x_5^3,
u_3=x_1^3x_2^3x_3^2x_4^3x_5^3,\\
&&u_4=x_1^3x_2^3x_3^3x_4^2x_5^3,
u_5=x_1^3x_2^3x_3^3x_4^3x_5^2.
\end{eqnarray*}
We use $c=x_1^3$ as a test element and see that, for example,
$cu_1^{7} \notin J^{[7]}$.  Similar
calculations and the observation in (\ref{notintc}) show that the
remaining monomials are not in $J^*$.  Computing $J:J^*$ gives the desired result.
\end{example}

\begin{remark}
Using our notation, the previous example is the case where $n=5$,
$d=4$ and $p=7$.
Since $d<n$, the results of Fedder and Watanabe \cite{FW}, Huneke
\cite{Hu2} and Hara \cite{Ha3} would predict that the
ring is $F$-regular if $p>11$, $p>8$ or $p>7$, respectively.
Note that in the previous example $p>d$, but $p$ is less than each of
the bounds and the ring is not $F$-regular.
\end{remark}

\begin{remark}
We have been able to compute the test ideal in the following cases:
$$
\begin{array}{l|l|l}
d&p&n \\ \hline
4 & 7 & 5 \\
5 & 2,3 & 3,4,5,6,7,8,9 \\
7 & 2,3,5 & 3,4,5,6,7,8,9 \\
8 & 3,5 & 3,4,5,6,7,8 \\
9 & 2 & 3,4,5,6,7,8,9 \\
9 & 5 & 3,4,5,6 \\
9 & 7 & 3,4 \\
10 & 3 & 3
\end{array}
$$
We have not included any examples where $p>d$ and we obtained the predicted
result, although we can compute many examples in those cases.  Our
computations tend to be limited by the degrees of the monomials
involved and the number of generators of the ideals involved.  As $p$
and $n$ grow, the degrees of the monomials grow, and as $n$ grows, the number
of generators of the ideals involved grows.  
Also, our examples do not represent the absolute limits of current
computation.  Computing further examples is incredibly time consuming
and for a fixed $p$ and $d$, the pattern as $n$ increases tends to
stabilize.  We expect that future results will eventually make further
computations unnecessary.

\end{remark}
\section{Non-integrally Closed Test Ideals}
Many of the test ideals that we can compute are not integrally
closed.  In this section we confirm that two of the examples of test
ideals in the previous section are not integrally closed.

\begin{example}
The test ideal computed in (\ref{n=d=5,p=2}) is not integrally closed.
Let $\tau$ be the test ideal.  The integral closure of $\tau$ is
$(x_1, \ldots, x_5)^3 \neq \tau$.  For example, $x_1x_2x_3 \in
(x_1, \ldots, x_5)^3 \backslash \tau$.
\end{example}

\begin{example}
The test ideal computed in (\ref{n=4,d=7,p=3}) is not integrally closed.
Let $\tau$ be the test ideal.  The integral closure of $\tau$ is
$(x_1, \ldots, x_4)^4 \neq \tau$.  For example, $x_1x_2^3 \in
(x_1, \ldots, x_4)^4 \backslash \tau$.
\end{example}

\bibliographystyle{alpha}

\end{document}